# Instant Evaluation and Demystification of $\zeta(n), L(n,\chi)$ that Euler, Ramanujan Missed-III

V.V.RANE

THE INSTITUTE OF SCIENCE,

15, MADAME CAMA ROAD,

MUMBAI-400 032,

INDIA

e-mail address : v_v_rane@yahoo.co.in

**Abstract** : For the Hurwitz zeta function $\zeta(s,\alpha)$, we show that for an integer $m \geq 0$, $\zeta(-m,\alpha)$, is a polynomial in $\alpha$ by using integration method unlike our earlier paper (with the same title), where we have proved the result by power series (in $\alpha$) method. However the expression for $\left(\frac{\partial}{\partial s}\zeta(s,\alpha)\right)_{s=-m}$ contains a term involving an integral, which cannot be computed explicitly. Consequently, Riemann zeta values at positive even integer arguments can be evaluated without involving Bernoulli numbers and those, at positive odd integer arguments cannot be computed explicitly. We also have corresponding results for Dirichlet L-series.

# Instant Evaluation and Demystification of $\zeta(n), L(n,\chi)$ that Euler, Ramanujan Missed-III

V.V.RANE

THE INSTITUTE OF SCIENCE,

15, MADAM CAMA ROAD,

MUMBAI-400 032,

INDIA

e-mail address : v_v_rane@yahoo.co.in

Let $s = \sigma + it$ be a complex variable, where $\sigma$ and t are real. For a complex number $\alpha \neq 0, -1, -2, \ldots\ldots$, let $\zeta(s,\alpha)$ be the Hurwitz zeta function defined by $\zeta(s,\alpha) = \sum_{n\geq 0}(n+\alpha)^{-s}$ for $\sigma > 1$ and its analytic continuation. For an integer $q \geq 1$, let $\chi(\mathrm{mod}\ q)$ be a Dirichlet character and let $L(s,\chi)$ be the corresponding Dirichlet L-series so that $L(s,\chi) = q^{-s}\sum_{a=1}^{q}\chi(a)\zeta(s,\frac{a}{q})$. In what follows, superscript 'dash' denotes differentiation with respect to s. Let $\tau(\chi) = \sum_{a=1}^{q}\chi(a)e^{\frac{2\pi i a}{q}}$. In author [1], we have proved the following Theorem 1 and Theorem 2 (independent of each other), which we state below as Proposition 1 and Proposition 2 respectively.

**Proposition 1**: We have for any integer $m \geq 0$ and for any complex variable $\alpha$,

$$\zeta(-m,\alpha) = \sum_{k=0}^{m}\binom{m}{k}\zeta(-k)\alpha^{m-k} + \alpha^m - \frac{\alpha^{m+1}}{m+1}, \text{ where } \binom{0}{0}=1.$$

**Proposition 2**: We have for any integer $m \geq 1$ and any complex $\alpha$,

$$\sum_{k=0}^{m-1}\binom{m}{k}\zeta(-k,\alpha) + \alpha^m = \frac{1}{m+1}.$$

Propositions 1 and 2 were obtained by considering the Taylor series expansion of the analytic function $\zeta(s,\alpha)$ (for a fixed s) as a function of $\alpha$ about 1 for $|\alpha-1|<1$; and the Taylor series of $\zeta(s,1+\alpha)$ about $\alpha$ for $|\alpha|<1$ respectively. In the case of Proposition 2, putting $m=1,2,3,\ldots$ therein successively, we get the polynomial expression in $\alpha$ for $\zeta(0,\alpha),\zeta(-1,\alpha),\zeta(-2,\alpha),\ldots$ respectively. However, we shall see that the expression for $\zeta'(-m,\alpha)$ for any integer $m\geq 0$, contains a term involving an integral, which cannot be computed and thus for $m\geq 0$, $\zeta'(-m,\alpha)$ cannot be computed explicitly and thus, following Lemma 0 below (which is Lemma 3 of author [1]), we cannot compute explicitly $\zeta(2n+1)$ for any integer $n\geq 1$; or we cannot compute $L(n,\chi)$, when an integer $n>1$ and $\chi\pmod{q}$ are of opposite parity.

We need following lemmas.

**Lemma 0 :** I) We have for $m\geq 1$, $\zeta(2m+1) = \frac{(-1)^m 2^{2m+1} \pi^{2m}}{(2m)!} \zeta'(-2m)$

II) If $\chi\pmod{q}$ is an even, primitive character, then we for $m\geq 0$.

$$L(2m+1,\chi) = \frac{(-1)^m 2^{2m+1} \pi^{2m} q^{-2m}}{\tau(\overline{\chi})(2m)!} L'(-2m,\overline{\chi}) = \frac{(-1)^m 2^{2m+1} \pi^{2m}}{\tau(\overline{\chi})(2m)!} \sum_{a=1}^{q} \overline{\chi}(a)\zeta'(-2m,\tfrac{a}{q})$$

III) If $\chi\pmod{q}$ is an odd, primitive character, then we have for $m\geq 1$,

$$L(2m,\chi) = \frac{i(-1)^{m+1} 2^{2m} \pi^{2m-1} q^{1-2m}}{\tau(\overline{\chi})(2m-1)!} L'(1-2m,\overline{\chi}) = \frac{i(-1)^{m+1} 2^{2m} \pi^{2m-1}}{\tau(\overline{\chi})(2m-1)!} \sum_{a=1}^{q} \overline{\chi}(a)\zeta'(1-2m,\tfrac{a}{q})$$

**Lemma 1 :** Let $\psi(u) = u-[u]-\tfrac{1}{2}$, where [u] is integral part of a real variable u. If $\sigma>0, x>0$ and $0<\alpha\leq 1$, then we have

I) $\zeta(s,\alpha) - \frac{x^{1-s}}{s-1} = \sum_{0 \leq n \leq x-\alpha} (n+\alpha)^{-s} + \psi(x-\alpha)x^{-s} - s\int_x^\infty \frac{\psi(u-\alpha)}{u^{s+1}} du.$

II) $\zeta'(s,\alpha) - \frac{d}{ds}\left(\frac{x^{1-s}}{s-1}\right) = \sum_{0 \leq n \leq x-\alpha} \frac{d}{ds}(n+\alpha)^{-s} + \psi(x-\alpha)\frac{d}{ds}x^{-s} + \frac{d}{ds}\int_x^\infty \psi(u-\alpha)\frac{\partial}{\partial u} u^{-s} du,$

where $\zeta'(s,\alpha) = \frac{\partial}{\partial s}\zeta(s,\alpha)$.

That is, $\zeta'(s,\alpha)$, the derivative of $\zeta(s,\alpha)$ with respect to s, can be obtained by term by differentiation (with respect to s) of the expression for $\zeta(s,\alpha)$ above in Lemma 1 (I). This result is valid for higher order derivatives (with respect to s) too. It is pertinent to note that in author [2], we have proved that the power series (in $\alpha$) of $\zeta'(s,1+\alpha)$ can be obtained by term-by-term differentiation (with respect to s) of the power series in $\alpha$ of $\zeta(s,1+\alpha)$.

Our Lemma 1 (I) partly illustrates the more general principle that if $F(s) = \sum_{n=1}^\infty f(n,s)$ is an analytic function of s in some right half plane of the complex plane, where $f(u,s)$ is a smooth function of u on an interval $(u_0, \infty)$ for some $u_0 \geq 0$, and of complex variable s, then the analytic continuation, functional equation or approximate function of F(s) if any, can be obtained by using just two facts namely,

1) Euler's summation formula

2) $\psi(u) = -\sum_{|n| \geq 1} \frac{e^{2\pi i n u}}{2\pi i n}$, the series converging boundedly on the interval $[0,1]$.

**Proof : I)** The proof follows by considering for $\sigma > 1$,

$\zeta(s,\alpha) = \sum_{n \geq 0}(n+\alpha)^{-s} = \left(\sum_{0 \leq n \leq x-\alpha} + \sum_{n > x-\alpha}\right)(n+\alpha)^{-s}$, where x>0 and where empty sum stands for zero. We then apply Euler's summation formula to the summation $\sum_{n>x-\alpha}(n+\alpha)^{-s}$.

For x>0 and for $0<\alpha\leq 1$, note that $x-\alpha>-1$.

By Euler's summation formula, we have for $\sigma>1$,

$$\sum_{n>x-\alpha}(n+\alpha)^{-s} = \int_{x-\alpha}^{\infty}(u+\alpha)^{-s}du - s\int_{x-\alpha}^{\infty}\psi(u)(u+\alpha)^{-s-1}du + x^{-s}\psi(x-\alpha)$$

$$=\int_x^{\infty} u^{-s}du - s\int_x^{\infty}\psi(u-\alpha)u^{-s-1}du + x^{-s}\psi(x-\alpha) = \frac{x^{1-s}}{s-1} - s\int_x^{\infty}\psi(u-\alpha)u^{-s-1}du + x^{-s}\psi(x-\alpha).$$

The proof of I) follows now.

II) We consider for $\sigma>1$ and for $x>0$,

$$\zeta'(s,\alpha) = -\sum_{n\geq 0}(n+\alpha)^{-s}\log(n+\alpha) = -\left(\sum_{0\leq n\leq x-\alpha} + \sum_{n>x-\alpha}\right)(n+\alpha)^{-s}\log(n+\alpha).$$

By Euler's summation formula we have for $\sigma>1$,

$$-\sum_{n>x-\alpha}(n+\alpha)^{-s}\log(n+\alpha) = -\int_{x-\alpha}^{\infty}(u+\alpha)^{-s}\log(u+\alpha)du$$

$$-\int_{x-\alpha}^{\infty}\psi(u)\tfrac{d}{du}(u+\alpha)^{-s}\log(u+\alpha))du - (x^{-s}\log x)\psi(x-\alpha)$$

$$= \int_{x-\alpha}^{\infty}\tfrac{\partial}{\partial s}(u+\alpha)^{-s}du + \int_{x-\alpha}^{\infty}\psi(u)\tfrac{\partial}{\partial u}\tfrac{\partial}{\partial s}(u+\alpha)^{-s}du + \psi(x-\alpha)\tfrac{\partial}{\partial s}x^{-s}$$

$$= \int_{x-\alpha}^{\infty}\tfrac{\partial}{\partial s}(u+\alpha)^{-s}du + \int_{x-\alpha}^{\infty}\psi(u)\tfrac{\partial}{\partial s}\tfrac{\partial}{\partial u}(u+\alpha)^{-s}du + \psi(x-\alpha)\tfrac{\partial}{\partial s}x^{-s}$$

$$= \tfrac{\partial}{\partial s}\int_{x-\alpha}^{\infty}(u+\alpha)^{-s}du + \tfrac{\partial}{\partial s}\int_{x-\alpha}^{\infty}\psi(u)\tfrac{\partial}{\partial u}(u+\alpha)^{-s}du + \psi(x-\alpha)\tfrac{\partial}{\partial s}x^{-s}$$

$$= \tfrac{\partial}{\partial s}\int_x^{\infty} u^{-s}du + \tfrac{\partial}{\partial s}\int_x^{\infty}\psi(u-\alpha)\tfrac{\partial}{\partial u}u^{-s}du + \psi(x-\alpha)\tfrac{\partial}{\partial s}x^{-s}$$

$$= \tfrac{\partial}{\partial s}\left(\tfrac{x^{1-s}}{s-1}\right) + \tfrac{\partial}{\partial s}\int_x^\infty \psi(u-\alpha)\tfrac{\partial}{\partial u}u^{-s}\,du + \psi(x-\alpha)\tfrac{\partial}{\partial s}x^{-s}$$

**Corollary of Lemma 1** : Choosing $x=1$ we have for $\sigma > 0$.

I) $\zeta(s,\alpha) = \alpha^{-s} + \tfrac{1}{s-1} + (\tfrac{1}{2}-\alpha) - s\int_1^\infty \tfrac{\psi(u-\alpha)}{u^{s+1}}\,du$ .

II) $\zeta'(s,\alpha) = -\alpha^{-s}\log\alpha - \tfrac{1}{(s-1)^2} - \int_1^\infty \tfrac{\psi(u-\alpha)}{u^{s+1}}\,du + s\int_1^\infty \tfrac{\psi(u-\alpha)}{u^{s+1}}\log u\,du$ .

Next, we state our Theorem. In what follows, $\zeta'(s,\alpha) \equiv \tfrac{\partial}{\partial s}\zeta(s,\alpha)$ and the empty sum will stand for zero .

**Theorem** : Let $r \geq 0$ and $N \geq 1$ be integers. For any integer $\ell \geq 1$ and for $0 < \alpha \leq 1$, let

$$\psi_\ell(\alpha) = -\sum_{|n|\geq 1}\frac{e^{2\pi i n\alpha}}{(2\pi i n)^\ell} = -\frac{\zeta(1-\ell,\alpha)}{(\ell-1)!} = \frac{B_\ell(\alpha)}{\ell!}$$, $B_\ell(\alpha)$ being Bernoulli polynomial of degree $\ell$ .

Then for $\sigma > -r$ and for $0 < \alpha \leq 1$, we have

I) $\zeta(s,\alpha) = \displaystyle\sum_{0\leq n\leq N-\alpha}(n+\alpha)^{-s} + \tfrac{N^{1-s}}{s-1} + \sum_{l=0}^{r}\tfrac{(-1)^l}{l!}\tfrac{s(s+1)\ldots(s+l-1)}{N^{s+l}}\zeta(-l,\alpha) - s(s+1)\ldots(s+r)\int_N^\infty \tfrac{\psi_{r+1}(u-\alpha)}{u^{s+r+1}}du$

II) $\zeta'(s,\alpha) = -\displaystyle\sum_{0\leq n\leq N-\alpha}(n+\alpha)^{-s}\log(n+\alpha) + \tfrac{d}{ds}\left(\tfrac{N^{1-s}}{s-1}\right)$

$+ \displaystyle\sum_{\ell=0}^{r}\tfrac{(-1)^l}{l!}\tfrac{\zeta(-\ell,\alpha)}{N^l}\tfrac{d}{ds}\left(\tfrac{s(s+1)\ldots(s+\ell-1)}{N^s}\right) - \int_N^\infty \psi_{r+1}(u-\alpha)\tfrac{d}{ds}\left(s(s+1)\ldots(s+r)u^{-s-r-1}\right)du$

$= -\displaystyle\sum_{0\leq n\leq N-1}(n+\alpha)^{-s}\log(n+\alpha) - \tfrac{N^{1-s}}{s-1}\left(\tfrac{1}{s-1}+\log N\right) + \sum_{\ell=0}^{r}\tfrac{(-1)^\ell}{\ell!}\tfrac{\zeta(-\ell,\alpha)}{N^{s+\ell}}\prod_{j=0}^{\ell-1}(s+j)\left(\sum_{k=0}^{\ell-1}\tfrac{1}{s+k} - \log N\right)$

$$-\prod_{j=0}^{r}(s+j)\left(\sum_{k=0}^{r}\tfrac{1}{s+k}\right)\int_{N}^{\infty}\psi_{r+1}(u-\alpha)u^{-s-r-1}du+\prod_{j=0}^{r}(s+j)\int_{N}^{\infty}\psi_{r+1}(u-\alpha)u^{-s-r-1}\log u\,du\,.$$

**Note :** The statement of Theorem is to be understood under the convention that empty sum stands for 0 and empty product stands for 1 . Thus, the expression for $\zeta'(s,\alpha)$ can be obtained by term-by-term differentiation (with respect to s) of the expression for $\zeta(s,\alpha)$ and by differentiating under integral sign , if necessary .

Replacing $r \geq 0$ by $m \geq 0$ and by letting $s \to -m$ from right, we have the following Corollary .

**Corollary :** We have for any integer $m \geq 0$,

1) $\zeta(-m,\alpha) = \sum_{0 \leq n \leq N-\alpha}(n+\alpha)^m - \tfrac{N^{m+1}}{m+1} + \sum_{\ell=0}^{m}\binom{m}{\ell}\zeta(-\ell,\alpha)N^{m-\ell}$

so that for $m \geq 1$, $\sum_{0 \leq n \leq N-\alpha}(n+\alpha)^m - \tfrac{N^{m+1}}{m+1} + \sum_{\ell=0}^{m-1}\binom{m}{\ell}\zeta(-\ell,\alpha)N^{m-\ell} = 0$.

This statement is the same as the one obtained by putting N=1 namely,

$\sum_{\ell=0}^{m-1}\binom{m}{\ell}\zeta(-\ell,\alpha)+\alpha^m - \tfrac{1}{m+1} = 0$ for $m \geq 1$, which is our Proposition 2 above .

2) We have $\zeta'(-m,\alpha) = -\sum_{0 \leq n \leq N-1}(n+\alpha)^m \log(n+\alpha) + \tfrac{N^{m+1}}{m+1}(\log N - \tfrac{1}{m+1})$

$$-\sum_{l=0}^{m}N^{m-\ell}\binom{m}{\ell}\zeta(-\ell,\alpha)\left(\log N + \sum_{j=m-\ell+1}^{m}\tfrac{1}{j}\right)+(-1)^{m+1}m!\int_{N}^{\infty}\tfrac{\psi_{m+1}(u-\alpha)}{u}du\,.$$

**Note** : Here $\int_N^\infty \frac{\psi_{m+1}(u-\alpha)}{u} du = \int_N^\infty \frac{d\psi_{m+2}(u-\alpha)}{u} du = -\frac{\psi_{m+2}(N-\alpha)}{N} + \int_N^\infty \frac{\psi_{m+2}(u-\alpha)}{u^2} du$.

**Proof** : By Lemma 1, for $\sigma > 0$ and for $0 < \alpha \le 1$ and on choosing $x = N \ge 1$ an integer,

we have $\zeta(s,\alpha) = \frac{N^{1-s}}{s-1} + \sum_{0 \le n \le N-\alpha}(n+\alpha)^{-s} + \psi(N-\alpha)N^{-s} - s\int_N^\infty \frac{\psi_1(u-\alpha)}{u^{s+1}} du$

where $\psi_\ell(u) = -\sum_{|n| \ge 1} \frac{e^{2\pi i n u}}{(2\pi i n)^\ell}$ for any integer $\ell \ge 1$. Note $\psi(u) = \psi_1(u)$.

Consider $-s\int_N^\infty \frac{\psi_1(u-\alpha)}{u^{s+1}} du = -s\int_N^\infty \frac{d\psi_2(u-\alpha)}{u^{s+1}} du$.

On integration by parts, we have for $\sigma > 1$,

$-s\int_N^\infty \frac{\psi_1(u-\alpha)}{u^{s+1}} du = -s\left\{ \left[\frac{\psi_2(u-\alpha)}{u^{s+1}}\right]_{u=N}^\infty + (s+1)\int_N^\infty \psi_2(u-\alpha) u^{-s-2} du \right\}$

$= \frac{s\psi_2(N-\alpha)}{N^{s+1}} - s(s+1)\int_N^\infty \psi_2(u-\alpha)^{-s-2} du$.

Continuing in this manner, on integration by parts, we have

$-s\int_N^\infty \frac{\psi_1(u-\alpha)}{u^{s+1}} du = \frac{s\psi_2(N-\alpha)}{N^{s+1}} + \frac{s(s+1)}{N^{s+2}}\psi_3(N-\alpha) + \ldots\ldots + \frac{s(s+1)\ldots\ldots(s+r-1)}{N^{s+r}}\psi_{r+1}(N-\alpha)$

$- s(s+1)\ldots\ldots(s+r)\int_N^\infty \psi_{r+1}(u-\alpha) u^{-s-r-1} du$.

Note that $\psi_\ell(N-\alpha) = \psi_\ell(-\alpha) = (-1)^\ell \psi_\ell(\alpha) = (-1)^{\ell-1} \frac{\zeta(1-\ell,\alpha)}{(\ell-1)!}$ for any integer $\ell \ge 1$.

This gives for $\sigma > -r$,

$\zeta(s,\alpha) = \sum_{0 \le n \le N-\alpha}(n+\alpha)^{-s} + \frac{N^{1-s}}{s-1} + \sum_{\ell=0}^r (-1)^\ell \frac{s(s+1)\ldots\ldots(s+\ell-1)}{l!.N^{s+l}} \zeta(-\ell,\alpha) - s(s+1)\ldots\ldots(s+r)\int_N^\infty \frac{\psi_{r+1}(u-\alpha)du}{u^{s+r+1}}$

Next, we choose N=1. Letting $s \to -r$ from right through real values and noting that the last term involving the integral vanishes and on replacing r by m, we get

$$\zeta(-m,\alpha) = \alpha^m - \tfrac{1}{m+1} + \sum_{\ell=0}^{m}\binom{m}{\ell}\zeta(-\ell,\alpha).$$

That is, $\alpha^m - \tfrac{1}{m+1} + \sum_{\ell=0}^{m-1}\binom{m}{\ell}\zeta(-\ell,\alpha) = 0$.

II) By Lemma 1 on replacing x by an integer $N \geq 1$, we have for $\sigma > 0$,

$$\zeta'(s,\alpha) = -\sum_{0 \leq n \leq N-1}(n+\alpha)^{-s}\log(n+\alpha) + \tfrac{d}{ds}\left(\tfrac{N^{1-s}}{s-1}\right) + \psi(N-\alpha)\tfrac{d}{ds}N^{-s} - \tfrac{d}{ds}\int_{N}^{\infty}s\psi(u-\alpha)u^{-s-1}du$$

$$= -\sum_{0 \leq n \leq N-1}(n+\alpha)^{-s}\log(n+\alpha) + \tfrac{d}{ds}\left(\tfrac{N^{1-s}}{s-1}\right)$$

$$+ \tfrac{d}{ds}\left\{\sum_{\ell=0}^{r}\tfrac{(-1)^{\ell}}{\ell!}\tfrac{s(s+1)\ldots(s+\ell-1)}{N^{s+\ell}}\zeta(-\ell,\alpha) - s(s+1)\ldots(s+r)\int_{N}^{\infty}\tfrac{\psi_{r+1}(u-\alpha)}{u^{s+r+1}}du\right\} \text{ for } \sigma > -r.$$

This results in the expression for $\zeta'(s,\alpha)$ as stated in statement II) of our Theorem.

**Proof of corollary** : 1) Consider the expression I) of $\zeta(s,\alpha)$ for $\sigma > -r$ and for $0 < \alpha \leq 1$. We let $s \to -r$ from right through real values. Note that the last term (having integral as a factor) varishes for $s = -r$, as this term has zero as a factor. Replacing r by m, we get the expression for $\zeta(-m,\alpha)$.

2) Consider the expression II) of $\zeta'(s,\alpha)$ for $\sigma > -r$ and for $0 < \alpha \leq 1$. Consider last two terms involving integrals in this expression, when $s \to -r$ from right through real values. Note that the last term vanishes for $s = -r$, as it has zero as a factor. As $s \to -r$, the last but one term except minus sign $= \lim_{s \to -r} \left( \prod_{j=0}^{r} (s+j) \right) \left( \sum_{k=0}^{r} \frac{1}{s+k} \right) \int_{N}^{\infty} \psi_{r+1}(u-\alpha) u^{-s-r-1} du$

$= \lim_{s \to -r} \left( \prod_{j=0}^{r} (s+j) \right) \left( \frac{1}{s+r} \right) \int_{N}^{\infty} \psi_{r+1}(u-r) u^{-s-r-1} du$

$= \left( \prod_{j=0}^{r-1} (-r+j) \right) \int_{N}^{\infty} \frac{\psi_{r+1}(u-\alpha)}{u} du = (-1)^r r! \int_{N}^{\infty} \frac{\psi_{r+1}(u-\alpha)}{u} du$.

Remaining first three terms of $\zeta'(s,\alpha)$ can be dealt with easily. Replacing r by m, we get the above expression for $\zeta'(-m,\alpha)$.

### References

1. V.V. Rane, Instant Evaluation and Demystification of $\zeta(n), L(n,\chi)$ that Euler, Ramanujan Missed-I ( arXiv.org website)

2. V.V. Rane, Instant Evaluation and Demystification of $\zeta(n), L(n,\chi)$ that Euler, Ramanujan Missed-II ( arXiv.org website) .